\documentclass[preprint,12pt]{elsarticle}

\usepackage{graphicx}

\usepackage{amssymb}
\usepackage{amsmath}

\journal{Applied Mathematics and Computation}

 \newcommand{\ds}[0] \displaystyle
 \newcommand{\OL}[1] {\mbox{O}\left(#1\right)}
 
 \newcommand{\dtot}[2] { \frac{d {#1} } {d {#2}} }
 
 \newcommand{\dpar}[2] { \frac{\partial {#1} } {\partial {#2}} }
 \newcommand{\ndpar}[3] { \frac{\partial^{#3} {#1} } {\partial #2 ^{#3} }}

\begin{document}

\begin{frontmatter}



\title{Removing trailing tails and delays induced by artificial dissipation in Pad\'e numerical schemes for stable compacton collisions}


\author{Julio Garral\'on}
\author{Francisco Rus}
\author{Francisco R. Villatoro\corref{autor}}

\address
  {Dept. Lenguajes y Ciencias de la Computaci\'on, \\
   Universidad de M\'alaga, 29071 M\'alaga, Spain}

\cortext[autor]{Corresponding author. Tel.: 34-951952388; fax:
 +34-951952542. {\it E-mail address}: villa@lcc.uma.es}

\begin{abstract}
The numerical simulation of colliding solitary waves with compact support arising from the Rosenau--Hyman $K(n,n)$ equation requires the addition of artificial dissipation for stability in the majority of methods. The price to pay is the appearance of trailing tails, amplitude damping, and delays as the solution evolves. These undesirable effects can be corrected by properly counterbalancing two sources of artificial dissipation; this procedure is designed by using the slow time evolution of the parameters of the solitary waves under the presence of the dissipation determined by means of adiabatic perturbation methods. The validity of the tail removal methodology is demonstrated on a Pad\'e numerical scheme. The tails are completely removed leaving only a small compact ripple at the original position of their front, and the numerical stability of the scheme under compacton collisions is preserved, as shown by extensive numerical experiments for several values of $n$.
\end{abstract}


\end{frontmatter}

\section{Introduction}

Rosenau and Hyman~\cite{RosenauHyman1993} came across solitary waves with compact support, therein referred to as \emph{compactons}, while studying the nonlinear dispersion in the formation of patterns in liquid drops. These traveling waves arose as solutions of a particular generalization of the well-known Korteweg--de Vries equation, commonly denoted as $K(m,n)$, that reads
\begin{equation}
 u_t + (u^m)_x + (u^n)_{xxx} = 0, \qquad m > 0, \quad 1 < n \le 3,
\label{eq:Kmn}
\end{equation}
where $u(x,t)$ is the wave amplitude as a function of the spatial variable $x$, and time $t$. Current analytical studies of this equation~\cite{OlverRosenau1996,Rosenau1996,LiOlver1997,LiOlver1998,Rosenau2000,GalaktionovSvirshchevskii2007,GalaktionovPohozaev2008,AmbroseWright2010}, cannot deal with the interaction between compactons. Hence, numerical methods, such as pseudospectral schemes~\cite{RosenauHyman1993,HymanRosenau1998},
finite elements~\cite{IsmailTaha1998,GarralonRusEtAl2006},
finite differences~\cite{IsmailTaha1998,DeFrutosEtAl1995},
Pad\'e approximants~\cite{RusVillatoro2007,MihailaEtAl2010},
modified equations~\cite{RusVillatoro2008}
or particle schemes based on the dispersive-velocity method~\cite{ChertockLevy2001}, must be used.

The majority of the numerical methods for compacton equations require the addition of artificial viscosity in order to cope with collisions, since otherwise instabilities appear which may cause a blow up in the solution. In pseudospectral schemes a hyperviscosity term is used where a second-order, linear dissipative term affects only the high-frequency modes of the solution by using a high-pass filter~\cite{HymanRosenau1998}. In both finite element and finite difference schemes a fourth-order, linear dissipative term is frequently incorporated~\cite{IsmailTaha1998,GarralonRusEtAl2006,DeFrutosEtAl1995,MihailaEtAl2010,RusVillatoro2008}. The addition of artificial viscosity in Eq.~\eqref{eq:Kmn} distorts the original compacton, generating small trailing tails, amplitude damping and velocity losses. The development of a method for the numerical stabilization of the compactons in collisions without the appearance of tails is an open problem to be addressed in this paper.

Numerically-induced phenomena can be studied and corrected by means of the method of modified equations~\cite{RusVillatoro2008,VillatoroRamos1999}. For the analysis of these phenomena, perturbation methods~\cite{KevorkianCole1996} can be applied to the analysis of effects introduced by local truncation errors as perturbations of the original evolution equation. These methods have been successfully applied to ordinary differential equations~\cite{VillatoroRamos1999,JunkYang2004} and nonlinear evolution equations~\cite{HermanKnickerbocker1993,MarchantSmyth1996,Marchant2002}. For compacton equations, adiabatic perturbation methods have been considered for the $K(n,n)$ equation in Garral\'on and Villatoro~\cite{GarralonVillatoro2012}, extending previous results for the $n=2$ case in Pikovsky and Rosenau~\cite{PikovskyRosenau2006}, and Rus and Villatoro~\cite{RusVillatoro2009}. These techniques, which generalize previous results for solitons~\cite{FernandezFroeschleEtAl1979,Lamb1980,Biswas2010}, can be used to analytically determine the trailing tails introduced by the artificial viscosity, opening the possibility of their removal by means of correction terms modifying the original evolution equation.

In this paper adiabatic perturbation methods are applied to the $K(n,n)$ equation and a new numerical technique for trailing tail removal is introduced. The validity of this procedure is checked by means of a Pad\'e numerical method. Section~\ref{sec:perturbations} is devoted to the development of the tail removal technique based on the adiabatic perturbation to the numerical scheme with artificial viscosity. In section~\ref{sec:numericalmethod} we briefly present the numerical method based on Pad\'e approximants used to solve the $K(n,n)$ equation. The results of extensive numerical experiments are presented in Section~\ref{sec:results} for several values of $n$. Final conclusions are summarized in Section~\ref{sec:conclusions}.

\section{Adiabatic perturbations}
\label{sec:perturbations}

The adiabatic perturbation method is applied to the $K(n,n)$ equation with both second- and fourth-order linear dissipations. Such a method determines the slow time evolution of the parameters of the compactly supported solitary waves by using that of the invariants under the dissipative perturbation.

Let us consider the perturbed $K(n,n)$ equation given by
\begin{equation}
 {u}_{t} + (u^n)_{x} + (u^n)_{xxx} = \varepsilon\,\mathcal{P}(u),
\label{eq:Knn-Pert}
\end{equation}
\noindent
where $1 < n \le 3$, the perturbation $\mathcal{P}(u)$ is a function of $u$ and its spatial and temporal derivatives, and $|\varepsilon| \ll 1$ is a small parameter. After multiplying Eq.~\eqref{eq:Knn-Pert} by $u^n$ and integrating in space, the only non-null term in its left-hand side is the first one, resulting in
\begin{equation}
\frac{d}{dt} \int^{\infty}_{-\infty} \frac{u^{n+1}}{n+1}\,dx =
             \varepsilon\,\int^{\infty}_{-\infty} u^{n}\,\mathcal{P}(u)\,dx,
\label{eq:Pert-un}
\end{equation}
whose left-hand side is the temporal derivative of the second invariant of the unperturbed $K(n,n)$, i.e., exactly nil for  Eq.~(\ref{eq:Knn-Pert}) with $\varepsilon=0$. However, perturbations such that the right-hand side of Eq.~\eqref{eq:Pert-un} is non-zero result in the adiabatic evolution of parameters of the compacton solution of the unperturbed equation. Introducing the slow time scale $\tau = \epsilon \, t$, the compacton solution of the $K(n,n)$ can be written as
\begin{equation}
u_c(x, t, \tau)=\left\{ \frac{2 \,n \, c(\tau)}{n + 1}
    \cos^2 \left( \frac{n-1}{2 \, n} \ (x - c(\tau)\,  t)
    \right)\right\}^{1/(n-1)},
\label{eq:Uctau}
\end{equation}
for $|x-c(\tau)\,t| \le n \,\pi / (n-1)$, and $u_c(x, t, \tau)=0$ otherwise. Inserting this ansatz into the perturbed equation~\eqref{eq:Pert-un} the slow time evolution of the velocity $c(\tau)$ of the perturbed compacton can be calculated; note that the amplitude of the perturbed compacton is uniquely determined by this velocity.

Let us determine the evolution of $c(\tau)$ for perturbation given by
\begin{equation}
\varepsilon\,\mathcal{P}(u) = \alpha_2\,\,u_{xx} - \alpha_4\,\,u_{xxxx},
\label{eq:Perts}
\end{equation}
where $|\alpha_2|,|\alpha_4|\ll 1$ are small parameters; this perturbation is dissipative if $\alpha_2>0$ and $\alpha_4>0$. The substitution of Eqs.~\eqref{eq:Perts} and~\eqref{eq:Uctau} into Eq.~\eqref{eq:Pert-un} yields an ordinary differential equation for $c(\tau)$ written as
\begin{equation}
 c'(\tau)= -\frac{(n-1)^2}{n\,(n+3)} \,\alpha_2\,\,c(\tau) -\frac{(n-1)^3\,((n-3)\,n-1)}{(n-5)\,n^3\,(n+3)}\, \alpha_4\,c(\tau).
 \label{eq:ode}
\end{equation}
The solution of this equation is
\begin{equation}
 c(\tau) = c(0)\,\exp\left(-\frac{(n-1)^2}{n\,(n+3)}\,\alpha_2\, \tau
  -\frac{(n-1)^3\,((n-3)\,n-1)} {(n-5)\,n^3\,(n+3)}\,\alpha_4\,\, \tau\right),
   \label{eq:ode:sol}
\end{equation}
showing that, under the perturbation, the velocity, as well as the amplitude, decays and produces lags compared with the unperturbed compacton.

The perturbation~\eqref{eq:Perts} introduces trailing tails in the perturbed compacton not accounted for by the ansatz~\eqref{eq:Uctau}, since the perturbed equation does not possess solitary waves solutions; hence the amplitude of the initial pulse slowly decays as described by Eq.~\eqref{eq:ode:sol}. This decay results in the formation of a plateau behind the compacton; this behaviour is similar to that of the solitons of the Korteweg--de Vries equation under the same kind of perturbation~\cite{Lamb1980}. Let us note that the adiabatic perturbation method can be used to estimate the shape of this tail, as shown in Refs.~\cite{PikovskyRosenau2006} and~\cite{RusVillatoro2009} for the $K(2,2)$ equation, and Ref.~\cite{GarralonVillatoro2012} for the $K(n,n)$ equation.

The main contribution of this paper is the introduction of a new procedure for removing the trailing tail by properly adjusting the values of the small parameters $\alpha_2$ and $\alpha_4$, suggested by the fact that their contributions in Eq.~\eqref{eq:Perts} have opposite signs. Fixing the value of $\alpha_2$ as a function of $\alpha_4$ and $n$ yields
\begin{equation}
  \alpha_2(n,\alpha_4) =   - \frac{(n-1)\,((n-3)\,n-1)}{(n-5)\,n^2}\,\alpha_4.
\label{eq:condition}
\end{equation}
For $1<n\le 3$, the value of $\alpha_2(n,\alpha_4)/\alpha_4$ is positive, increasing from $1/9$ to $1/4$ as $n$ decreases from $3$ to $2$, and decreasing from $1/4$ to $0$ as $n$ decreases from $2$ to $1$.

\section{The numerical method}
\label{sec:numericalmethod}

The tail removal procedure presented in the previous section could be incorporated into a numerical method for Eq.~\eqref{eq:Kmn} which uses artificial viscosity in order to deal with compactons collisions. Obviously, this procedure could introduce new instabilities. For illustration purposes, let us consider one of the most used methods for Eq.~\eqref{eq:Kmn}, a Pad\'e approximation in space, with periodic
boundary conditions in a finite interval $[0,L]$, and a method of lines in time as described in~\cite{RusVillatoro2007,MihailaEtAl2010}. Equation~\eqref{eq:Knn-Pert} with perturbations~\eqref{eq:Perts} and a moving frame of reference with velocity $c_0$ is given by
\begin{equation}
   Eq[u] \equiv u_t - c_0\,u_x + (u^n)_x + (u^n)_{xxx} - \alpha_2\,u_{xx} + \alpha_4\,u_{xxxx} = 0.
\label{eq:movingKnn}
\end{equation}
The Pad\'e method of this equation results in
\begin{eqnarray}
 &&
        \mathcal{A}(\operatorname{E}) \dtot{U_j}{t}
   -c_0 \mathcal{B}(\operatorname{E}) U_j
   + \mathcal{B}(\operatorname{E}) (U_j)^n
   + \mathcal{C}(\operatorname{E}) (U_j)^n
 \nonumber \\ && \phantom{ \mathcal{A}(\operatorname{E}) \dtot{U_j}{t}}
   - \alpha_2\,\mathcal{S}(\operatorname{E}) U_j
   + \alpha_4\,\mathcal{D}(\operatorname{E}) U_j = 0, \qquad j=0, 1, \ldots M,
\label{eq:Knn-Pade24}
\end{eqnarray}
where the spatial grid nodes are $x_j = j\,\Delta x$, for $j=0, 1, \ldots M$, with $\Delta x=L/M$, the numerical solution is $U_j(t) \approx u(x_j, t)$, the shift operator $\operatorname{E}$ is defined as $\operatorname{E}\,U_j=U_{j+1}$, and $\mathcal{A}^{-1}(\operatorname{E})\,\mathcal{B}(\operatorname{E})$, $\mathcal{A}^{-1}(\operatorname{E})\,\mathcal{S}(\operatorname{E})$ $\mathcal{A}^{-1}(\operatorname{E})\,\mathcal{C}(\operatorname{E})$, and $\mathcal{A}^{-1}(\operatorname{E})\,\mathcal{D}(\operatorname{E})$ are Pad\'e operators for the first-, second-, third-, and fourth-order spatial derivatives respectively. Concretely, the Pad\'e operators to be used in this paper are
\begin{displaymath}
  \mathcal{A}(\operatorname{E}) = \frac{
              \operatorname{E}^{-2} + 26 \operatorname{E}^{-1} +
              66 +
              26\operatorname{E}^1 + \operatorname{E}^2} {120},
\label{eq:PadeApproxA}
\end{displaymath}
\begin{displaymath}
  \mathcal{B}(\operatorname{E}) = \frac{
             -\operatorname{E}^{-2} -10 \operatorname{E}^{-1}
             +10\operatorname{E}^1 +\operatorname{E}^2} {24\Delta x},
\label{eq:PadeApproxB}
\end{displaymath}
\begin{displaymath}
  \mathcal{S}(\operatorname{E}) = \frac{
              \operatorname{E}^{-2} +2\operatorname{E}^{-1}
             -6
             +2\operatorname{E}^1 +\operatorname{E}^2} {6\Delta x^2},
\label{eq:PadeApproxS}
\end{displaymath}
\begin{displaymath}
  \mathcal{C}(\operatorname{E}) = \frac{
             -\operatorname{E}^{-2} + 2\operatorname{E}^{-1}
             -2\operatorname{E}^1 + \operatorname{E}^2} {2\Delta x^3},
\label{eq:PadeApproxC}
\end{displaymath}
and
\begin{displaymath}
  \mathcal{D}(\operatorname{E}) = \frac{
              \operatorname{E}^{-2} -4\operatorname{E}^{-1}
             +6
             -4\operatorname{E}^1 +\operatorname{E}^2} {\Delta x^4},
\label{eq:PadeApproxD}
\end{displaymath}
corresponding to the following approximations to the first four spatial derivatives
\begin{displaymath}
  \mathcal{A}^{-1}(\operatorname{E})\,\mathcal{B}(\operatorname{E})\,u(x_j,t) =
   \dpar{u}{x}(x_j,t) + \frac{\Delta x^6}{5040}\,\ndpar{u}{x}{7}(x_j,t) + \OL{\Delta x^8},
\label{eq:PadeApproxAB}
\end{displaymath}
\begin{displaymath}
  \mathcal{A}^{-1}(\operatorname{E})\,\mathcal{S}(\operatorname{E})\,u(x_j,t) =
   \ndpar{u}{x}{2}(x_j,t) + \frac{\Delta x^4}{720}\,\ndpar{u}{x}{6}(x_j,t) + \OL{\Delta x^6},
\label{eq:PadeApproxAS}
\end{displaymath}
\begin{displaymath}
  \mathcal{A}^{-1}(\operatorname{E})\,\mathcal{C}(\operatorname{E})\,u(x_j,t) =
   \ndpar{u}{x}{3}(x_j,t) - \frac{\Delta x^4}{240}\,\ndpar{u}{x}{7}(x_j,t) + \OL{\Delta x^6},
\label{eq:PadeApproxAS}
\end{displaymath}
and
\begin{displaymath}
  \mathcal{A}^{-1}(\operatorname{E})\,\mathcal{D}(\operatorname{E})\,u(x_j,t) =
   \ndpar{u}{x}{4}(x_j,t) - \frac{\Delta x^2}{12}\,\ndpar{u}{x}{6}(x_j,t) + \OL{\Delta x^4}.
\label{eq:PadeApproxAS}
\end{displaymath}
The modified equation of the numerical method~\eqref{eq:Knn-Pade24}, which includes the local truncation error terms, is given by
\begin{equation}
   Eq[u(x_j,t)]- \alpha_4\,\frac{\Delta x^2}{12}\,\ndpar{u}{x}{6}(x_j,t) + \OL{\Delta x^4}  = 0,
\label{eq:movingKnn:modified}
\end{equation}
showing that the numerical scheme is dissipative (for $\alpha_4>0$) after the application of the tail removal procedure.

For the discretization in time of Eq.~\eqref{eq:Knn-Pade24}, let us use the second-order accurate, implicit midpoint rule, yielding
\begin{eqnarray}
  &&
  \mathcal{A}(\operatorname{E}) \frac{U_j^{(i+1)}-U_j^{(i)}}{\Delta t} - \left(
  c_0 \mathcal{B}(\operatorname{E})
    + \alpha_2\,\mathcal{S}(\operatorname{E})
    - \alpha_4\,\mathcal{D}(\operatorname{E}) \right)\,
   \frac{U_j^{(i+1)}+U_j^{(i)}}{2}
  \nonumber \\ &&
  \phantom{  \mathcal{A}(\operatorname{E}) \frac{U_j^{(i+1)}-U_j^{(i)}}{\Delta t} }
  +(\mathcal{B}(\operatorname{E}) + \mathcal{C}(\operatorname{E}))
     \left( \frac{U_j^{(i+1)}+U_j^{(i)}}{2}\right)^n  =  0,
\label{eq:Knn-midpoint}
\end{eqnarray}
where $\Delta t$ is the time step and $U_j^{(i)} \approx u(x_j, i\,\Delta t)$. The resulting nonlinear system of equations is solved for $U_j^{(i+1)}$ by using the Newton's iterative method.

\section{Presentation of results}
\label{sec:results}

Let us study the behaviour of the new procedure for trailing tails removal on the compacton solutions of the $K(n,n)$ equation with $n \in \{3$, 2, 5/3, 3/2, 7/5, 4/3, 9/7, $5/4 \}$. Numerical solutions with several values of $\alpha_4$ will be compared with both $\alpha_2=0$ and $\alpha_2=\alpha_2(n,\alpha_4)$, for the propagation of a one-compacton solution and for the collisions between two compactons.

\subsection{One-compacton solutions}

\begin{figure}
\centering
{\includegraphics[width=14cm]{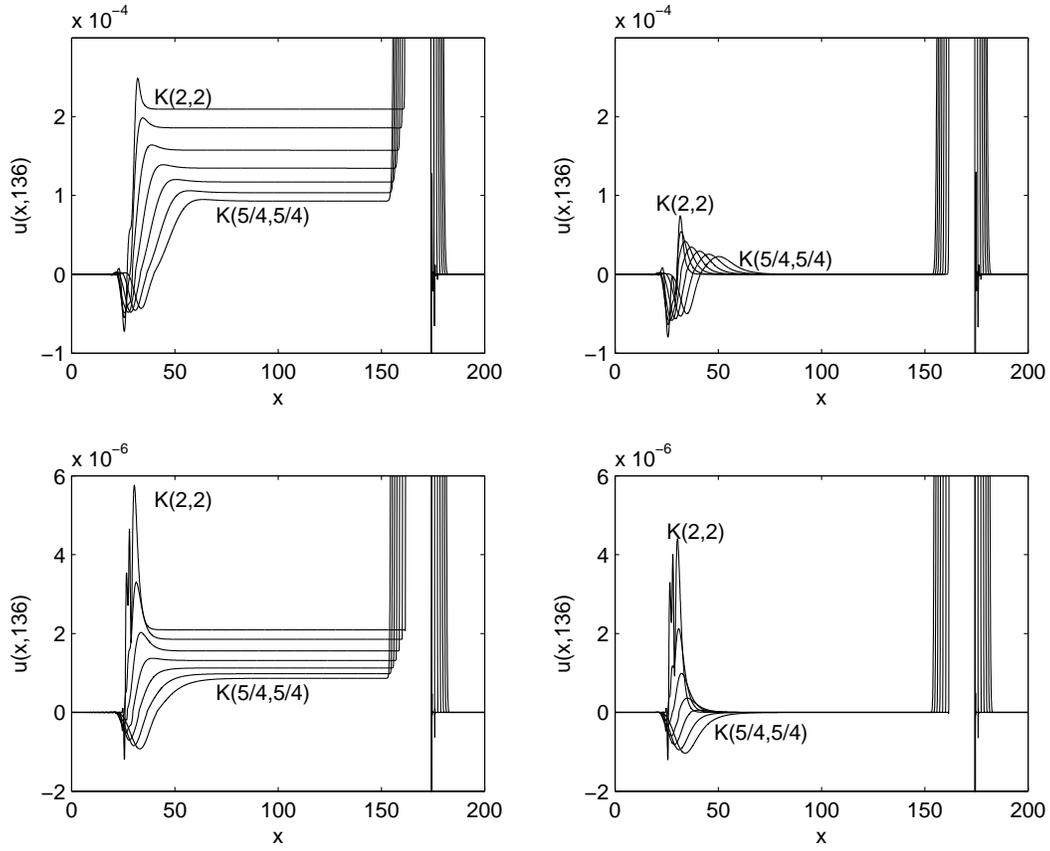} }
\caption{Numerical one-compacton solution of the $K(n,n)$ equation for $n = 2$, 5/3, 3/2, 7/5, 4/3, 9/7, and $5/4$, with $c=1$, $c_0=0.5$, $\Delta x = 0.05$, and $\Delta t = 0.04$. The top plots use $\alpha_4=10^{-3}$ and the bottom ones $\alpha_4=10^{-5}$; the left plots use $\alpha_2=0$ and the right ones $\alpha_2(n,\alpha_4)$.}
\label{fig:colas}
\end{figure}

Let us first consider the evolution of a one-compacton solution with and without the tail removal procedure. Figure~\ref{fig:colas} shows snapshots of the numerical solution at the same instant of time for several values of $n$. The left plots show the trailing tails that appear when $\alpha_2=0$, for $\alpha_4=10^{-3}$ (top plot) and $\alpha_4=10^{-5}$ (bottom plot). Note that a zoom in has been used since the amplitude of the tails is several orders of magnitude smaller than that of the compacton. The tail starts at the initial location of the compacton at $t=0$ with a ripple with a negative peak followed by a positive one; after the ripple, the tail has a nearly constant plateau connecting with the left edge of the compacton. The amplitude of the tail decreases as $n$ does and as so does $\alpha_4$. The right plots show the removal of these tails when using $\alpha_2(n,\alpha_4)$, for $\alpha_4=10^{-3}$ (top plot) and $\alpha_4=10^{-5}$ (bottom plot). The plateau completely disappears, but the ripple in the tail front persists with smaller amplitude. In long-time integrations, the compacton is far away and clearly separated from the ripple located at its initial position, so it behaves as if there is no tail.

\begin{figure}
\centering
{ \includegraphics[width=7cm]{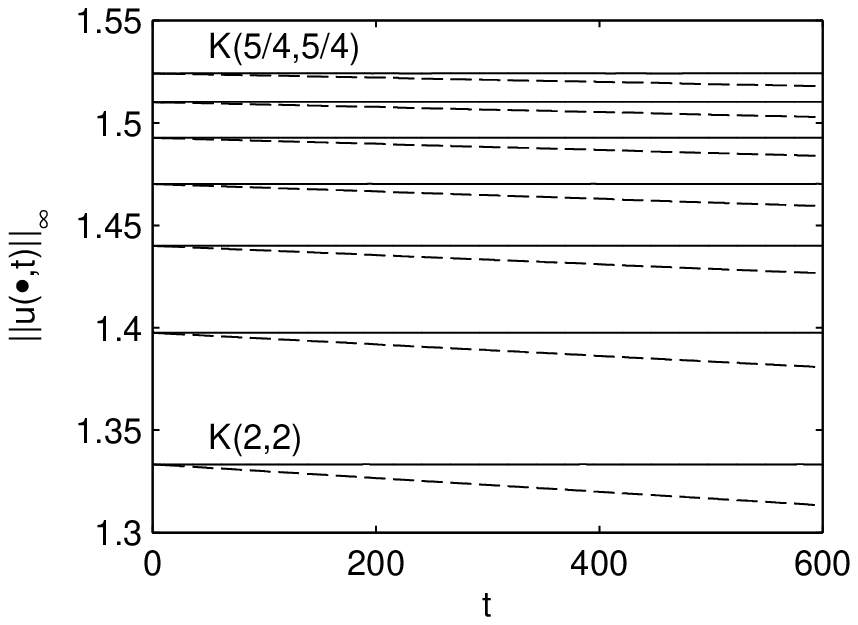} }
\caption{The evolution in time of the maximum amplitude of the one-compacton solution of the $K(n,n)$ equation for
$n = 2$, 5/3, 3/2, 7/5, 4/3, 9/7, and $5/4$, with $c=1$, $\alpha_4=10^{-3}$, $c_0=1$, $\Delta x = 0.05$, and $\Delta t = 0.04$. The dashed lines use $\alpha_2=0$ and the solid ones $\alpha_2(n,\alpha_4)$.}
\label{fig:amplitudes}
\end{figure}

Figure~\ref{fig:amplitudes} compares the evolution in time of the compacton maximum amplitude with $\alpha_2=0$ (dashed lines) and $\alpha_2(n,\alpha_4)$ (solid ones), i.e., with and without the tails. The appearance of the tails compensates the loss in the compacton maximum amplitude due to the dissipation since the first invariant is preserved during the propagation, hence, with the tail removal procedure, the maximum amplitude practically retains its initial value.

Table~\ref{tab:delays} shows the difference between the location of the analytical and numerical compacton maximum amplitudes at a time $t=2000$ for several values of $\alpha_4$, with $\alpha_2=0$ and $\alpha_2(n,\alpha_4)$. The numerical velocity of the compactons in the dissipativeless case (when $\alpha_2=\alpha_4=0$) is very near, but not exactly equal, to the analytical one, resulting in a very small difference in the position at $t=2000$ of about $0.6$ (for the values of $\Delta x$ and $\Delta t$ used in Table~\ref{tab:delays}), which is mostly independent of the value of $n$. The use of the tail removal procedure, $\alpha_2(n,\alpha_4$), results in exactly the same value for this difference in location at $t=2000$, except for either $n=3$, or $\alpha_4=10^{-2}$ and $n\ge 2$, as shown in Table~\ref{tab:delays}. Without the tail removal procedure, the difference between the numerical and analytical compactons after 20000 time steps can be very large, in fact, larger than 26, 3, and 0.8 (in units of space) for $\alpha_4=10^{-2}$, $10^{-3}$, and $10^{-4}$ respectively. Hence, Table~\ref{tab:delays} shows that the tail removal procedure works properly and results in a numerical velocity equal to that of the numerical method without dissipation.

\begin{table}
\begin{center}
\begin{tabular}{c|cc|cc|cc}
\hline
&
    \multicolumn{2}{c|}{$\alpha_4=10^{-2}$} &
    \multicolumn{2}{c|}{$\alpha_4=10^{-3}$} &
    \multicolumn{2}{c}{$\alpha_4=10^{-4}$} \\
\cline{2-7}
 $n$ & $\alpha_2=0$ & $\alpha_2(n,\alpha_4)$ &
  $\alpha_2=0$ & $\alpha_2(n,\alpha_4)$ &
  $\alpha_2=0$ & $\alpha_2(n,\alpha_4)$ \\
\hline
   3 & 505.50 & 108.50 & 53.20 & 5.50 & 5.80 & 0.90 \\
   2 & 404.90 & 0.90 & 47.80 & 0.60 & 5.40 & 0.60 \\
 5/3 & 223.10 & 0.70 & 25.20 & 0.60 & 3.10 & 0.60 \\
 3/2 & 129.80 & 0.60 & 14.40 & 0.60 & 2.00 & 0.60 \\
 7/5 & 80.80 & 0.60 & 9.00 & 0.60 & 1.40 & 0.60 \\
 4/3 & 53.20 & 0.60 & 6.00 & 0.60 & 1.10 & 0.60 \\
 9/7 & 36.70 & 0.60 & 4.20 & 0.60 & 0.90 & 0.60 \\
 5/4 & 26.30 & 0.60 & 3.20 & 0.60 & 0.80 & 0.60 \\
\hline
\end{tabular}
\caption{Difference between the maximum amplitude location of the analytical and numerical one-compacton solutions at time $t=2000$ for $\alpha_4=10^{-2}$, $10^{-3}$, and $10^{-4}$, with $\alpha_2=0$ and $\alpha_2(n,\alpha_4)$. In this table, $c=1.0$, $c_0=0.5$, $\Delta x = 0.1$, and $\Delta t = 0.1$ have been used.}
\label{tab:delays}
\end{center}
\end{table}

\subsection{Compacton collisions}

The numerical simulation of compacton collisions requires the use of artificial viscosity in order to avoid the appearance of instabilities. It seems that the tail removal procedure introduced in this paper may  affect the stability of the simulations for colliding compactons since the explicit artificial viscosity terms has been canceled. Let us summarize our extensive set of simulations showing this is not the case.

\begin{figure}
\centering
{ \includegraphics[width=14cm]{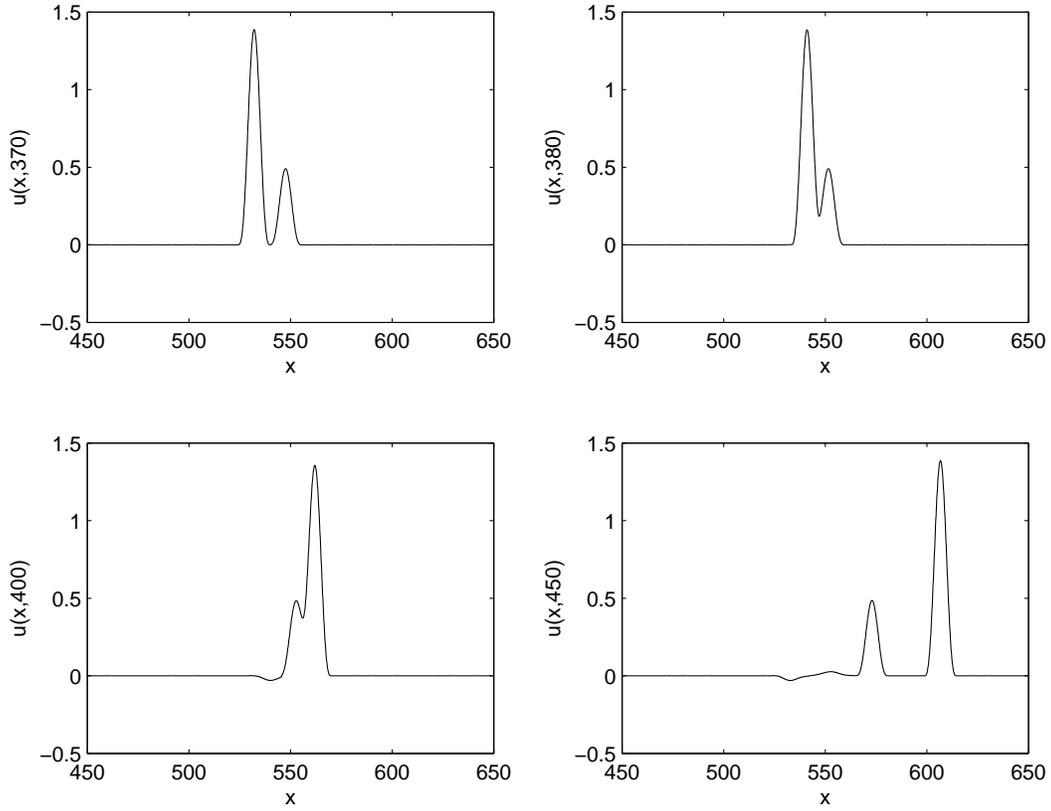} }
\caption{Four snapshots of the collision of two compactons of the $K(5/3,5/3)$ equation with velocities $c_1=1.0$ and $c_2=0.5$ at times $t=370$, 380, 400, and 450, from upper left to lower right. The parameters used in this simulation are $\alpha_2=0$, $\alpha_4=10^{-3}$, $c_0=0.1$, $\Delta x = 0.1$, and $\Delta t = 0.1$.}
\label{fig:collision}
\end{figure}

Figure~\ref{fig:collision} illustrates the collision of two compactons of the $K(5/3,5/3)$ equation with velocities $c_1=1.0$, the tallest one, and $c_2=0.5$. The upper left plot shows both compactons before the collision. Two snapshots of the collision are presented in the upper right and the lower left plots. Finally, the lower right plot indicates that both compactons recover their amplitude and velocities, but a small ripple remains at the location of the interaction. This collision is typical for all the $K(n,n)$ equations.

\begin{figure}
\centering
{ \includegraphics[width=14cm]{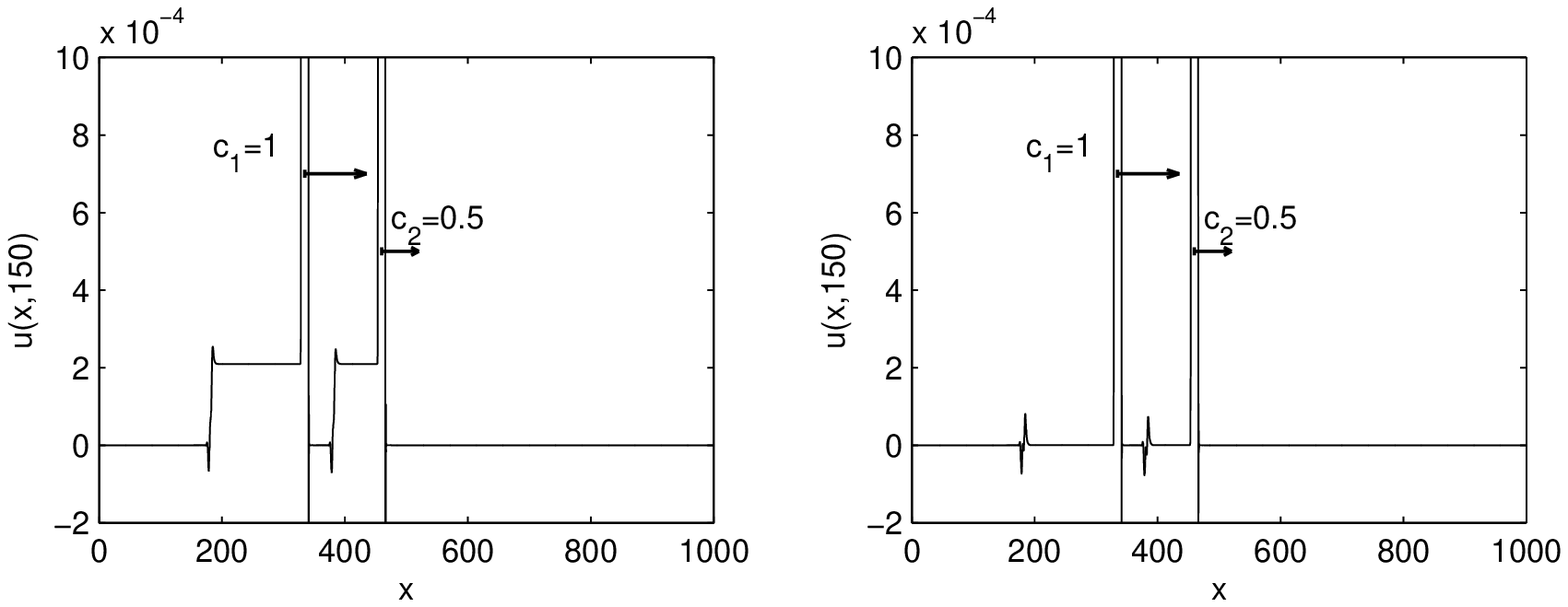} }
{ \includegraphics[width=14cm]{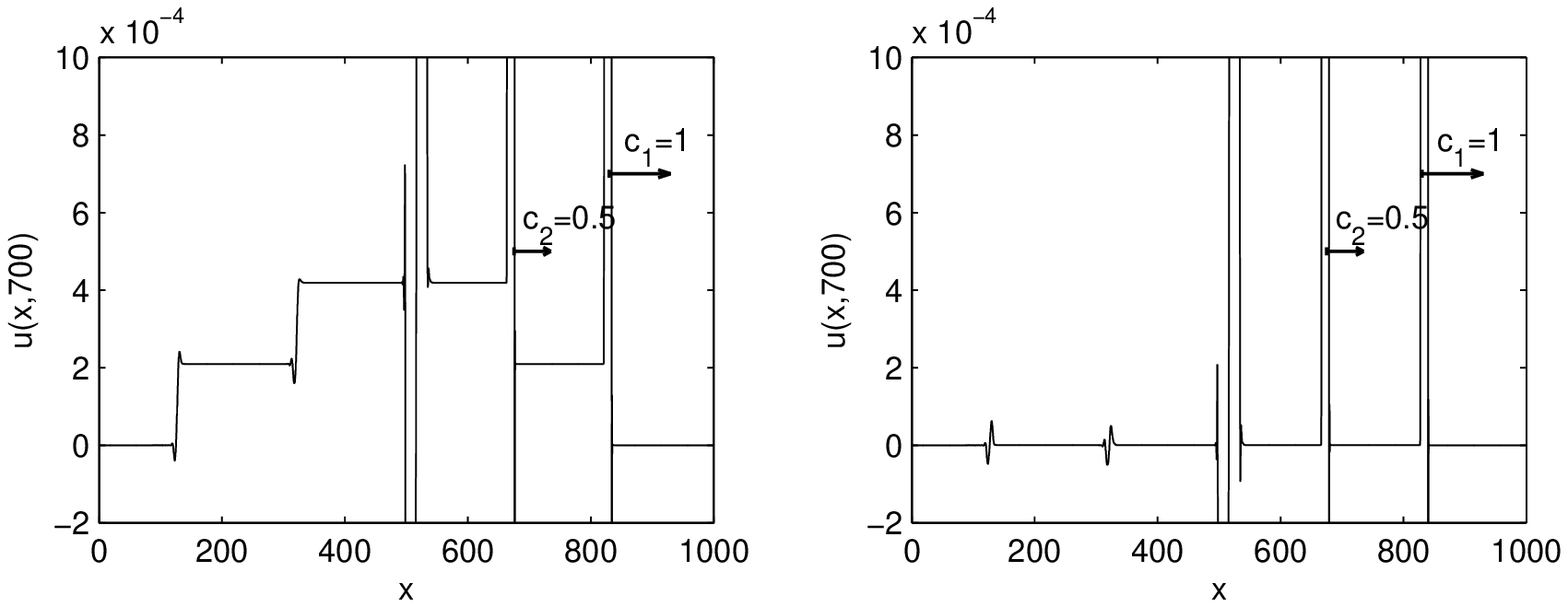} }
{ \includegraphics[width=14cm]{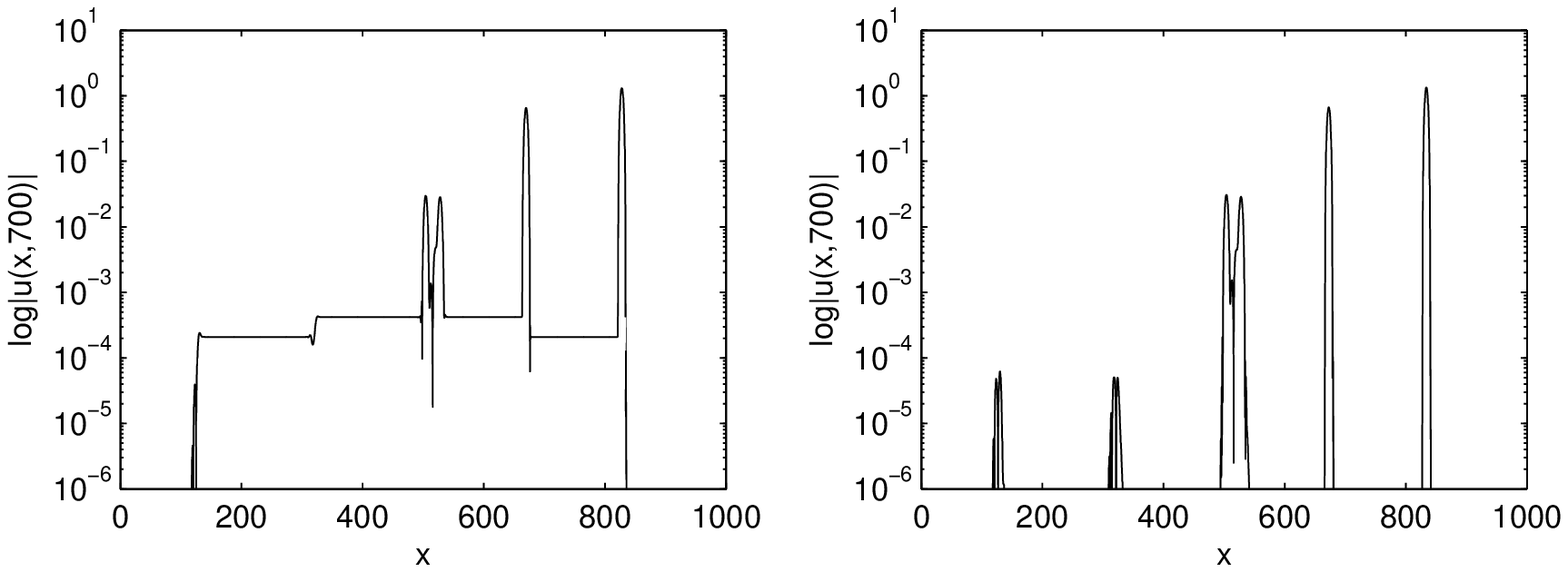} }
\caption{Numerical simulation of two colliding compactons of the $K(2,2)$ equation with velocities $c_1=1.0$ and $c_2=0.5$ before (top plots with $t=150$) and after (middle and bottom ones with $t=700$) their mutual collision, with $\alpha_2=0$ (left plots) and $\alpha_2(n,\alpha_4)$ (right ones). The parameters used in this simulation are $\alpha_4=10^{-3}$, $c_0=0.1$, $\Delta x = 0.1$, and $\Delta t = 0.1$.}
\label{fig:collisiontails}
\end{figure}

Figure~\ref{fig:collisiontails} shows a zoom in of the numerical simulation of the collision between two compactons of the $K(2,2)$ equation. The initial positions of the compactons and their velocities ($c_1=1$, $c_2=0.5$ and $c_0=0.1$) have been selected in order to ensure that the fronts of the tail of each compacton are far away at the time of collision. Left plots show both compactons for $\alpha_0=0$, before (top plot) and after (middle and bottom ones) their collision. The amplitude of the tail of each compacton depends only on the value of $\alpha_4$ (equal to $10^{-3}$ in this figure) when the tails are well separated (see upper left plot in Figure~\ref{fig:collisiontails}), but when the fastest compacton gets to the trailing tail of the slowest, the tail of this second one rides over that of the first one (see middle left plot). As it is widely known in compacton literature~\cite{RosenauHyman1993,GarralonRusEtAl2006,DeFrutosEtAl1995,CardenasEtAl2011}, a ripple  appears after the collision which can be seen around $x=500$ in the lower left plot; the ripple amplitude is about two orders of magnitude larger than that of the tails.

The application of the tail removal procedure cancels the trailing tails without affecting the stability of the numerical method during the collision of the compactons, as illustrated in the right plots shown in Fig.~\ref{fig:collisiontails} for the $K(2,2)$ equation. The small ripples located at the front of the tails which is not removed by taking $\alpha_2(n,\alpha_4)$ (see upper right plot) do not affect the collision which occurs as if they were absent. Lower right plot illustrates that the compactons, the ripple associated to the collision, and those associated to the tails behave as compactly supported solutions (except for the numerically induced, self-similar, backward and forward radiations reported by Rus and Villatoro~\cite{RusVillatoro2007a}, whose amplitude is smaller than $10^{-6}$ in the plots of Fig.~\ref{fig:collisiontails}).

\begin{figure}
\centering
{ \includegraphics[width=14cm]{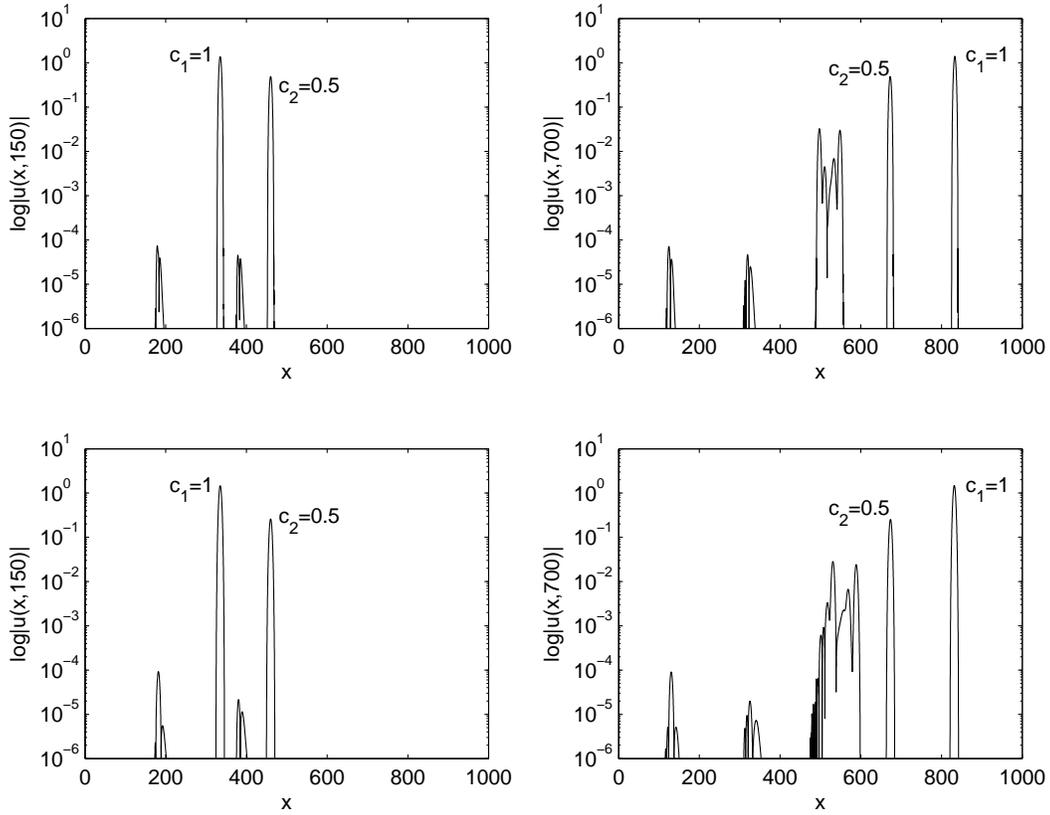} }
\caption{Numerical simulation of two colliding compactons of the $K(5/3,5/3)$ (top plots) and $K(7/5,7/5)$ (bottom ones) equations with velocities $c_1=1.0$ and $c_2=0.5$ before (left plots) and after (right ones) their mutual collision, with $\alpha_2(n,\alpha_4)$. The parameters used in this simulation are $\alpha_4=10^{-3}$, $c_0=0.1$, $\Delta x = 0.1$, and $\Delta t = 0.1$.}
\label{fig:collisiontailsbis}
\end{figure}

The results shown in Figure~\ref{fig:collisiontails} for the $K(2,2)$ equation are representative of an extensive set of simulations of compactons collisions for $n \in \{3$, 2, 5/3, 3/2, 7/5, 4/3, 9/7, $5/4 \}$, different values of $\alpha_4$, $c_1$, $c_2$, $c_0$, $\Delta x$, and $\Delta t$. Figure~\ref{fig:collisiontailsbis} shows representative results for the $K(5/3,5/3)$ (top plots) and $K(7/5,7/5)$ (bottom ones) equations. The slowest compacton has been stopped by using $c_0$ and the plots show that the tails have been removed, except for the ripple at their fronts, which are two orders of magnitude smaller than the residual after the collision of the compactons. Further results for other values of $n$ are omitted here for the sake of brevity. In all the cases the proper behavior of the tail removal procedure has been observed.

In long-time integrations, using periodic boundary conditions, two compactons collide multiple times and their trailing tails ride over one another increasing their total amplitude. In such simulations, the value of the artificial viscosity $\alpha_4$ must be chosen properly in order to avoid the appearance of instabilities resulting in the sudden blow-up of the solution inside the time integration interval. Our results show that the tail removal procedure do not affect the time of blow-up, validating its good performance in terms of the preservation of the stability of the numerical scheme.

\section{Conclusions}
\label{sec:conclusions}

A procedure for trailing tail removal for numerical methods for the $K(n,n)$ equation incorporating artificial viscosity has been introduced. This procedure is based on the analysis of the effect of the artificial viscosity in the propagation of compactons by means of the adiabatic perturbation method. The performance of the new procedure is illustrated by using a widely used numerical method based on Pad\'e approximants (which also can be derived from finite element and finite difference formulations). The results obtained after an extensive set of simulations show the effectiveness of new the trailing tail removal in the propagation of solutions with both only one compacton and two compactons in mutual interaction.

The new procedure could be applied to other numerical schemes for the $K(n,n)$ equation and for the numerical study of other nonlinear evolution equations presenting compactly supported solutions, such as those with cosine/sine compactons reported in Ref.~\cite{RusVillatoro2009b}.

\section*{Acknowledgements}

The authors would like to thank the anonymous reviewers for their valuable comments and suggestions to improve the quality of the paper. The research reported here has been partially supported by Projects
MTM2010--19969, TIN2008--05941 of the Ministerio de Ciencia e
Innovaci\'on of Spain, and TIC-6083 of the Junta de Andaluc\'ia of Spain.


\newcommand{\bookref}[5]{ {#4}, {#1}, {#2}, #3.}

\newcommand{\paperref}[7]{ {#1}, {#2}, {#7} {#5} (#4) {#6}.}
\newcommand{\paperrefno}[8]{#1, {{#3}} {\bf #4}  (#7) p. #8.}

\end{document}